\documentclass {article}
\usepackage{authblk}
\usepackage{blindtext}

\usepackage{latexsym,amsfonts,amsmath,amssymb,amsthm}

\usepackage{graphicx}
\usepackage{wrapfig}
\usepackage{bm}          
\usepackage[font=small,labelfont=bf,margin=10pt]{caption}
\usepackage[usenames,dvipsnames]{color} 
\usepackage{floatrow}

\definecolor{LinkColor} {rgb}{0, 0, .8} 
\usepackage[colorlinks=true, linkcolor=LinkColor, citecolor=black, urlcolor=blue,
    bookmarks, bookmarksnumbered=false, bookmarksopen, bookmarksopenlevel=2,
    verbose=false,
    pdfpagemode=UseOutlines,
    naturalnames,
    hypertexnames,
    pagebackref,
    hyperindex,
    plainpages=false,
    pdfpagelabels,
    pdftitle={Sobczyk},
    pdfauthor={Alan Macdonald}] {hyperref}

\usepackage{shadethm}
\definecolor{shadethmcolor}{gray}{.9}

\newcommand     {\separate}{\medskip}
\newcommand     {\half} {\frac{1}{2}}

\renewcommand {\outer}   {\mathbin{\wedge}}

\renewcommand {\i}  {\mathbf{i}}
\renewcommand {\j}  {\mathbf{j}}
\renewcommand {\k} {\mathbf{k}}
\renewcommand {\r} {\mathbf{r}}
\newcommand     {\s} {\mathbf{s}}

\usepackage[textwidth=12cm,textheight=22.8cm]{geometry}

\begin{document}
\date{}

\title{Sobczyk's simplicial calculus does\\ not have a proper foundation}
\author[1]{Alan Macdonald}
\affil[1]{Luther College, Decorah, IA 52101, USA}
\affil[ ]{macdonal@luther.edu}
\maketitle
\vspace{-.6in}
\renewcommand{\abstractname}{}
\begin{abstract}
\noindent
{The pseudoscalars in Garret Sobczyk's paper \emph{Simplicial Calculus with Geometric Algebra} \cite{Sobczyk} are not well defined.
Therefore his calculus does not have a proper foundation.}

\vspace{.02in}
Update: Sobczyk has conceded the point.
See the end of this paper.
\end{abstract}

\pagenumbering{gobble}

\separate
In 1868 Joseph Serret defined the area of a surface as the limit of areas of inscribed polyhedra \cite{serret}.
It is a natural generalization of the definition of the length of a curve as the limit of lengths of inscribed polygons.
Unfortunately, Serret's definition is flawed:
In 1882 Hermann Schwarz \cite{schwarz} showed that Serret's definition does not give a well defined surface area.
His simple example is the area of the lateral side of a cylinder.%
\footnote{One correct geometric definition of surface area divides a surface into small pieces, projects each piece onto the tangent plane at some point in the piece, adds the areas of the projections, and takes the limit as the pieces become smaller \cite{Franklin}.
The usual area formula $\iint_A |\r_u \times \r_v|\,dA$, where  $\r\colon A \to S$ is a regular parameterization of a surface,
is then a theorem.
This approach is too advanced for an introductory vector calculus course,
so the integral is usually taken as the definition of surface area.}

This story is told in an award winning paper by Frieda Zames \cite{Zames}
and illustrated in a Wolfram Demonstrations Project \cite{Wolfram}.
I learned of the problem with Serret's definition as an undergraduate in David Widder's advanced calculus text \cite{Widder},
and was quite intrigued.

\separate
Every manifold has a pseudoscalar field $I(x)$.
Garret Sobczyk's paper \emph{Simplicial Calculus with Geometric Algebra} attempts a geometric definition of $I(x)$ \cite{Sobczyk}.
It has been cited approvingly by Hestenes \cite{Hestenes} and used implicitly by Doran and Lasenby \cite{Doran}.
However, Sobczyk's definition is flawed for the same reason that Serret's definition is flawed: It is not well defined.

Sobczyk's definition of $I(x)$ is one of ``the formal definitions upon which we construct our theory''.
Since the definition is not sound, neither is the theory.
Following Schwarz, I will show that Sobczyk's definition does not give a well defined $I(x)$
for the lateral side of a cylinder.

\begin{wrapfigure}{r}{1.3in}
\caption{Triangles\\ inscribed in cylinder.}
\label{fig:Zames}
{\includegraphics[width=1.3in]{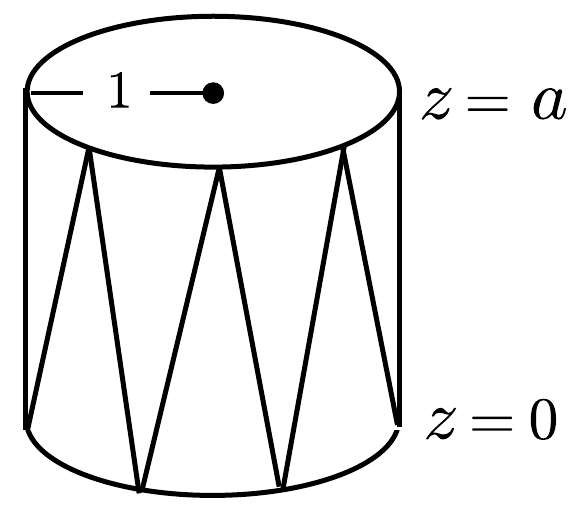}}
\end{wrapfigure}


Start with a cylinder of height $a$ erected on the circle $x^2 + y^2 = 1$ in the $z = 0$ plane.
See Fig.~\ref{fig:Zames}.
Inscribe a polyhedron in the cylinder whose faces consist of isosceles triangles as follows.
Divide the circumference at the top and bottom into $n$ equal arcs subtending angles $2\pi/n$ from the center,
but let the points of subdivision at the top lie midway between those at the bottom.

Draw a straight line from each point to its two neighbors on the same circle and to the
two nearest points of subdivision on the other circle.
The inscribed polyhedron thus formed has $2n$ isosceles triangles for faces.
\clearpage

\begin{wrapfigure}{r}{1.6in}
\vspace{-.15in}
\centering \includegraphics[width=1.6in]{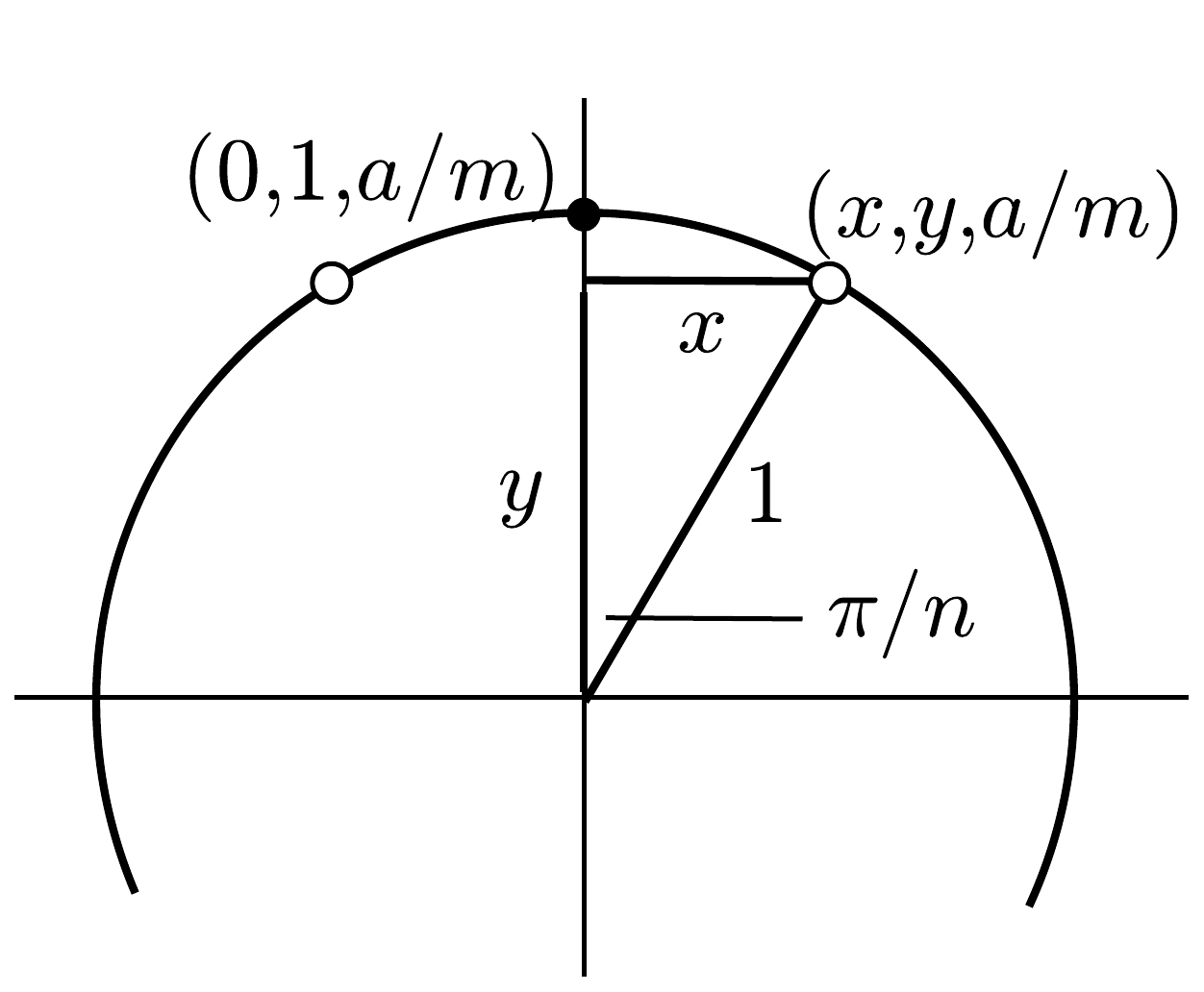}
\caption{The $z = a/m$ plane.}
\label{fig:SobczykFig}
\end{wrapfigure}
Next, suppose that the number of sides of the polyhedron is increased
by first dividing the cylinder into $m$ equal cylinders by planes parallel
to the base and then proceeding with each as above.

Consider, for example, an ``upside down'' isosceles triangle with apex in the $z = 0$ plane at $(0,1,0)$ and base in the $z = a/m$ plane.
Figure~\ref{fig:SobczykFig} shows that plane.
The two vertices on the base, denoted $\circ$, have coordinates
\[ \Big(\!\!\pm x,\, y,\, \frac{a}{m}\Big) = \Big(\!\!\pm\sin\frac{\pi}{n},\,  \cos\frac{\pi}{n},\,  \frac{a}{m}\Big). \]

Let $\r$ and $\s$ be the vectors from the apex $(0, 1, 0)$ to the base angles.
The plane and area of the triangle are represented by $\half \r \outer \s$.
Sobczyk \cite{Sobczyk} defines the pseudoscalar at $(0, 1, 0)$ to be
    \[ \lim_{\half\r \outer \s \to 0}  \frac{\half \r \outer \s}{|\half \r \outer \s|}
                            = \lim_{\r \outer \s \to 0}  \frac{\r \outer \s}{|\r \outer \s|}. \]

I now show that the limit is not well defined: it depends on the way one takes $\r \outer \s \to 0$.
Compute:
\begin{align*}
\r \outer \s
  &= \Big\{\!\sin\!\frac{\pi}{n}\,\i + \Big(\!\cos\!\frac{\pi}{n} - 1\Big)\j + \frac{a}{m}\k\Big\}
     \!\outer\! \Big\{\!\!-\sin\!\frac{\pi}{n}\,\i   + \Big(\!\cos\!\frac{\pi}{n} - 1\Big)\j + \frac{a}{m}\k\Big\} \\
  &= 2\sin\!\frac{\pi}{n}\,\Big(\!\cos\!\frac{\pi}{n} - 1\Big)\,\i \outer \j
     \,+\, \frac{2a}{m}\sin\!\frac{\pi}{n}\,\i \outer \k \\
  &\approx\frac{\pi}{n}\Big(\frac{\pi}{n}\Big)^2 \j \outer \i
     \,+\, \frac{2a}{m}\frac{\pi}{n}\, \i \outer \k \\
  &= \frac{\pi}{n}\Big\{\!\Big(\frac{\pi}{n}\Big)^2 \j \outer \i \,+\, \frac{2a}{m}\, \i \outer \k \Big\}.
\end{align*}

If $m = n$, then the $\i \outer \k$ term dominates for large $n$,
so $\lim_{n \to \infty}\frac{\r \outer \s}{|\r \outer \s|} = \i \outer \k$.
This is a unit pseudoscalar of the tangent algebra to the cylinder at the apex $(0,1,0)$.
Good! 

But if $m = n^3$, then the limit is $\j \outer \i$.
This is orthogonal to the pseudoscalar $\i \outer \k$ at $(0,1,0)$.
And if $m = cn^2$ ($c > 0$, a constant) then the limit is
a normalized linear combination of $\i \outer \k$ and $\j \outer \i$ with positive coefficients.

This does not yet contradict Sobczyk's claim,
as he requires that the vertices of the inscribed polyhedra for a given $n$ be contained in those for $n+1$.
For this, choose $n = 2^i$ and let $i \to \infty$.
And for the $m = cn^2$ case, use only powers of $2$ for $c$.

\separate
Of course the lateral side of a cylinder has a pseudoscalar field.
Unfortunately, Sobczyk's definition \cite{Sobczyk} does not capture it.

\separate
\textbf{Acknowledgements.}
I thank Professor Michael D. Taylor for helpful discussions.
\clearpage

Update: Sobczyk has conceded the point:
\href{https://gaupdate.wordpress.com/2017/10/24/a-macdonald-sobczyks-simplicial-calculus-does-not-have-a-proper-foundation/#comments}{Click}.
(The URL is long.)

\separate
Sobczyk's comment:
\begin{quote}
The flaw in my partitioning scheme had already been pointed out in
F. James,  ``Surface Area and the Cylinder Area Paradox'', The College Mathematics Journal 8 (1977) 207-211.
It is obviously not a trivial problem and it is one that is worth looking into.
If \mbox{MacDonald} [sic] had done that, he would have made a real contribution.
\end{quote}

\separate
My reply:
\begin{quote}Zames' paper (which I cite and describe) has nothing to do with geometric algebra. It describes a similarly failed construction using simplices, that of 2D (surface) area. Google Scholar lists 36 citations to Sobczyk's paper. As far as I know, the flaw in it had not previously been noted,
even by citing Zames.

There are correct geometric definitions of 2D area using simplices.
I cite one in a footnote.
As far as I know, none generalize to nD volume. So I don't see much hope for the more difficult task of salvaging Sobczyk's nD paper
while retaining ``simplicial'' in its title.
\end{quote}


\begin{thebibliography}{}
\small

\bibitem{Sobczyk}{G. Sobczyk, \emph{Simplicial calculus with Geometric Algebra}, Eq. (3.4).
In \emph{Clifford Algebras and their Applications in Mathematical Physics}, Springer (1992) 279-292.
(Also available at \url{http://geocalc.clas.asu.edu/pdf/SIMP_CAL.pdf})}

\bibitem{Hestenes}{D. Hestenes, \emph{Differential Forms in Geometric Calculus}, in: F. Brackx et al (eds.),
Clifford Algebras and their Applications in Mathematical Physics, Kluwer: \mbox{Dordercht}/Boston(1993), 269-285.}

\bibitem{Doran}{C. Doran and A. Lasenby, \emph{Geometric Algebra for Physicists},  Cambridge University Press, Cambridge, UK (2003), 203.}

\bibitem{serret}{J. A. Serret, \emph{Cours de Calcul Diff\'{e}rentiel et Integral (2 vol)}, Gauthier-Villars, Paris (1868) 296.}

\bibitem{schwarz}{H. A. Schwarz, \emph{Sur une definition  erron\'ee de l'aire, d'une surface courbe, Ges. Math. Abhandl. I}, 309.}

\bibitem{Zames}{F. Zames, \emph{Surface Area and the Cylinder Area Paradox}, The College Mathematics Journal \textbf{8} (1977) 207-211.
(Also avaiable as ``surface area.pdf'' at\\ \url{https://www.math.washington.edu/~morrow/334_13}.)}

\bibitem{Wolfram}{\url{http://demonstrations.wolfram.com/CylinderAreaParadox}.}

\bibitem{Widder}{D. Widder, \emph{Advanced Calculus, 2nd Ed.}, Dover (1989) 204.}

\bibitem{Franklin}{P. Franklin, \emph{A Treatice on Advanced Calculus}, Dover (2016) 373-376.}

\end{thebibliography}
\end{document}